\documentclass[11pt]{article}

\usepackage{amssymb,amsbsy,amsmath,amscd,amsthm,amsfonts}
\usepackage{cite}
\usepackage[utf8]{inputenc}
\usepackage{enumerate,hyperref}
\usepackage{dsfont}

\hypersetup{
    colorlinks=true,       
}

\newtheorem{theorem}{Theorem}
\newtheorem{condition}{Assumption}

\newtheorem{lemma}{Lemma}
\newtheorem{remark}{Remark}

\newtheorem{corollary}{Corollary}

\newcommand{\R}{\mathbb R}
\newcommand{\PP}{\mathbb P}
\newcommand{\E}{\mathbb E}
\newcommand{\VV}{\mathbb V}

\title{\textbf{Convex Set Detection}}
\author{Victor-Emmanuel Brunel\\
   CREST, Paris (France) - University of Haifa, Haifa (Israel)\\
   \texttt{victor.emmanuel.brunel@ensae-paristech.fr}}

\date{}

\begin{document}

\maketitle

\begin{abstract}
    We address the problem of one dimensional segment detection and estimation, in a regression setup. At each point of a fixed or random design, one observes whether that point belongs to the unknown segment or not, up to some additional noise. We try to understand what the minimal size of the segment is so it can be accurately seen by some statistical procedure, and how this minimal size depends on some \textit{a priori} knowledge about the location of the unknown segment. 
\end{abstract}


\section{Introduction}
	
	Consider the statistical model: 
	\begin{equation}
		\label{Model}
		Y_i=\mathds 1(X_i\in G)+\xi_i, i=1,\ldots,n,
	\end{equation}
	where $\mathds 1(\cdot)$ is the indicator function.
	The set made of the points $X_i$ is called the design.
	The unknown set $G$ is a segment on $[0,1]$, the noise terms $\xi_i$ are i.i.d. random variables, independent of the design. We distinguish two types of design:
	\begin{enumerate}
		\item[(DD)] Deterministic, and regular design: $X_i=i/n, i=1,\ldots,n$;
		\item[(RD)] Random, uniform design: the variables $X_i, i=1,\ldots,n$, are i.i.d. uniform on $[0,1]$.
	\end{enumerate}
	In the sequel, the design will be denoted by $\mathcal X$.
	This model can be interpreted as a partial and noisy observation of an image, $[0,1]$, in which there is an unknown object $G$. From this observation, one would like to determine whether it is true that there is an object - the unknown set $G$ might be empty -, and/or to recover that object, i.e. to estimate $G$. Another framework for this model is the noisy observation of some signal, here $\mathds 1(X_i\in G)$, and one would like to determine if the observation comes from a pure noise, or if there actually is some signal, and/or to recover that signal.
	
	As it was mentioned in \cite{Brunel2013}, what makes difficult estimation of $G$ is the complexity of the class of possible candidates, and detectability of $G$. The unknown set belongs to the class of all segments $[a,b]$ of $[0,1]$. This is a parametric class, with two parameters $0\leq a\leq b\leq 1$. Therefore, one may expect to be able to estimate $G$ at the parametric speed $1/n$, up to a positive multiplicative constant. However, it turns out that $G$ may be too hard to detect and thus, the speed of estimation would be deteriorated. For instance, if the segment is too small, and if no point of the design falls inside $G$, no procedure would see $G$. We would like to understand under which assumptions detection is not an obstacle for estimating $G$. In particular, the two following assumptions will be of interest for us :
	
\begin{condition}\label{A1}
	The set $G$ is of the form $[0,\theta]$, for some unknown number $0\leq \theta \leq 1$.
\end{condition}

\begin{condition}\label{A2}
	$|G|\geq \mu$. 
\end{condition}	
Here and in the rest of the paper, $\mu\in(0,1)$ is a given positive number. The brackets $|\cdot|$ stand for the Lebesgue measure. 

Let us discuss these two assumptions in order to understand, intuitively, why they facilitate detection of $G$. The first one gives some information on the location of $G$. In other words, it tells that the set $G$ starts from the left side of the frame. The second one tells that $G$ is not too small, and thus should not be unnoticed by the statistician. Model \eqref{Model}, together with Assumption \ref{A1}, is well known under the name of the change point problem. It can be rewritten as :
\begin{equation}
	\label{ChangePoint}
	Y_i=\mathds 1(X_i\leq \theta)+\xi_i, i=1,\ldots,n,
\end{equation}
for some number $\theta\in [0,1]$. The change point problem was studied in \cite[Sec. 1.9]{KTlectureNotes1993}, and a continuous-time version of this model is addressed in \cite{Korostelev1986}. The aim is to estimate the breakpoint $\theta$. 
In the continuous time version, Korostelev \cite{Korostelev1986} proposed a more general framework. Instead of the indicator function in the regression equation \eqref{ChangePoint}, that is a function with a jump at the point $\theta$, and satisfying a Lipschitz condition both on the left and on the right sides of $\theta$. In these two works \cite{Korostelev1986},\cite{KTlectureNotes1993}, the change point $\theta$ is estimated with a precision, in expectation, of the order of $1/n$, if $\theta$ is assumed to be separated from 0 and 1 : $h\leq \theta \leq 1-h$, for some $h\in (0,1/2)$. Ibragimov and Khasminskii \cite{Ibragimov1984}, under, among others, the same assumption of separation from 0 and 1, also proposed a consistent estimator of the discontinuity point of a regression function, with precision of order $1/n$ as well. This separation hypothesis is common in this kind of estimation problems. For instance, it is made in \cite{KorostelevTsybakov1992} and \cite[Chap 3]{KTlectureNotes1993}, where Tsybakov and Korostelev propose an estimator for boundary fragments. A boundary fragment, in dimension $d\geq 1$, is a subset $G$ of $\R^d$ which can be described as the subgraph of a positive function $g:[0,1]^{d-1} \longrightarrow [0,1]$:
\begin{equation}
	\label{BF}
	G=\left\{(x_1,\ldots,x_d)\in[0,1]^d : 0\leq x_d\leq g(x_1,\ldots,x_{d-1})\right\}.
\end{equation}
The authors study a similar model to Model \eqref{Model}, where $G$ is a boundary fragment instead of a segment. This is one possible generalization in higher dimensions of our problem. They assume that the true $g$ belongs to some Hölder class, and is separated from 0 and 1:
\begin{equation*}
	h\leq g(x)\leq 1-h, \forall x\in[0,1]^{d-1},
\end{equation*}
for a given parameter $h\in (0,1/2)$. In \cite{KorostelevTsybakov1992}, the same authors estimate the support of a uniform density, assuming it is a boundary fragment. Again, they make the hypothesis that the underlying function $g$ is separated from 0 and 1. In both models, they build a piecewise polynomial estimator of the function $g$, and prove that it is optimal in a minimax sense - see more details bellow -. One goal of the present work is to understand whether this separation hypothesis is necessary, in the simple case of dimension 1.

Two other extensions in higher dimension have been addressed in \cite{Brunel2013}. A segment of $[0,1]$ can be interpreted both as a one dimensional convex body, or as a one dimensional convex polytope. If the dimension is greater than 1, then the class of convex bodies is much bigger than that of convex polytopes. This is what explains that the optimal speed of estimation is better for polytopes than for general convex bodies - details in \cite{Brunel2013} -. In this previous work, we did not study the impact of such assumptions as Assumptions \ref{A1} and \ref{A2}. This question is beyond the scope of the present paper, but we believe that the results would be similar in the case of convex polytopes. Indeed, that case allows to keep the parametric property and, as we mentioned in \cite[Section 5]{Brunel2013}, detectability becomes the main obstacle for estimation of the unknown set. 
In \cite{Gayraud2001}, a testing problem is addressed, in a different model: the set of interest, as in \cite{KorostelevTsybakov1992}, is the support of a uniform density, assumed to be a boundary fragment (see \eqref{BF}). Of course, detection is not pertinent in that model, since, as long as there are observations, there must be a nonempty support. However, the author studied the separation rate for distinguishing hypotheses of the type $H_0 : G=G_0$ and $H_1 : d(G,G_0)\geq h$, where $G_0$ is a given boundary fragment, $d$ is the Nykodim distance between sets, defined below, and $h$ is a given positive number which may depend on $n$.

In \cite{ChanWalther2013}, the detection question is addressed in a slightly different framework. In Model \eqref{Model}, we assume that the strength of the signal is given, equal to 1. It is also interesting to deal with the case of a signal of unknown strength, i.e.
	\begin{equation}
		\label{Model2}
		Y_i=\delta\mathds 1(X_i\in G)+\xi_i, i=1,\ldots,n,
	\end{equation}
where $\delta$ is a positive number. For the signal to be detectable, there should be a tradeoff between its length $|G|$ and its strength $\delta$. Naturally, if $\delta$ is small, then the set $G$ should be big enough and conversely, if $\delta$ is large, the set $G$ is allowed to be small, so that the signal can be detected. Testing the presence of a signal, i.e. whether $\delta=0$ or not, is considered in \cite{ChanWalther2013}. This work is concerned with the power of two tests: the scan - or maximum - likelihood ratio, and the average likelihood ratio. The two tests are compared in two regimes: signals of small scales, i.e. $|G|\longrightarrow 0$, and signals of large scales, i.e. $\displaystyle{\underset{n\rightarrow\infty}{\operatorname{liminf}} \text{ } |G|>0}$. The design is (DD), and it is proved that if $\displaystyle{\delta\sqrt{n|G|}\geq\sqrt{2\ln\frac{1}{|G|}}+b_n}$, for some sequence $b_n$ such that $b_n\longrightarrow \infty$, then there is a test whose power is asymptotically 1. Note that here, $\delta\sqrt{|G|}$ is exactly the $L^2$-norm of the signal, $\|\delta\mathds 1(\cdot\in G)\|_2=\delta\sqrt{|G|}$. In \cite{LepskiTsybakov2000}, signals of unknown shape but known smoothness were considered. Exact minimax separation rates, in terms of the $L^2$-norm of the signal, for distinguishing the null hypothesis, under which observations are pure noise, and the alternative one, under which there is a signal, are given. Detection is harder in that framework, because unlike in Model \eqref{Model} or \eqref{Model2}, where the shape of the signal is known - it is piecewise constant -, only its smoothness is known, and the separation rates are larger than those of models \eqref{Model} and \eqref{Model2}, in the sense that they allow less freedom for the size of the signal. However, this problem is different from ours, since we are concerned with the location of the signal, not the signal itself. 

Model \ref{Model} deals with change points in the mean of the observations, conditionally to the design. Under Assumption \ref{A1}, there is only one change point, and under Assumption \ref{A2}, there are two. A problem of interest in time series analysis is that of detecting change points in the mean, or in the covariance function, or other characteristics of the series. We refer to \cite{ShaoZhang2010} and the references therein. In \cite{FrickMunkSieling2013}, a sample of $n$ independent observations $Y_1,\ldots,Y_n$ is given, and one assumes that $Y_i$ admits a density $f(\cdot,\vartheta(\frac{i}{n}))$ with respect to a given measure, for $i=1,\ldots,n$, where $f$ belong to an exponential parametric class of densities. The real valued function $\vartheta$ is assumed to be piecewise constant on $[0,1]$, with a finite number $K$ of jumps, not necessarily known. Under a similar condition on $\delta$ and $|G|$ to that of \cite{ChanWalther2013}, it is shown that at least one change point is consistently detectable: the proposed estimator of $K$ is positive with probability that goes to 1, as $n\rightarrow\infty$. When $f(\cdot,\theta)$ is the density of a Gaussian distribution with mean $\theta$ and given variance $\sigma^2>0$, the problem was addressed in \cite{Lebarbier2003}. This model includes models \eqref{Model} and \eqref{Model2}, when the design is (DD). However, Assumption \ref{A2} is not considered in that work, and optimality of the estimator of the change point is not treated, in the case when it is assumed to be unique - which corresponds to Assumption \ref{A1} - or when it is known that there are only two of them.



\subsection{Notation}

	If $G$ is a segment of $[0,1]$, we denote, respectively, by $\PP_G$, $\E_G$ and $\VV_G$ the joint probability measure of the observations $\left((X_1,Y_1),\ldots,(X_n,Y_n)\right)$ that come from Model \eqref{Model}, and the corresponding expectation and variance operators. If the event inside the probability, expectation or variance sign does not depend on $G$, we omit the subscript and write only $\PP$, $\E$ and $\VV$. 
	
	If $G_1$ and $G_2$ are two segments of $[0,1]$, we denote by $G_1\triangle G_2$ their symmetric difference. The Lebesgue measure of the symmetric difference of $G_1$ and $G_2$, $|G_1\triangle G_2|$, is also called the Nykodim distance between $G_1$ and $G_2$.

In the present work, we do two types of inference on the underlying set $G$: hypothesis testing and estimation. The tests consist of deciding whether $G$ is empty or not. Estimation consists of giving an estimator of $G$. A subset $\hat G_n$ of $[0,1]$, whose construction depends on the sample is called a set estimator or, more simply, an estimator.
Given an estimator $\hat G_n$, we measure its accuracy in a minimax framework. The risk of $\hat G_n$ on a class $\mathcal C$ of Borel subsets of $\R$ is defined as
\begin{equation*}
	\mathcal R_n(\hat G_n ; \mathcal C) = \sup_{G\in\mathcal C}\mathbb E_G[|G\triangle\hat G_n|].\tag{$*$}
\end{equation*}
The rate (a sequence depending on $n$) of an estimator on a class $\mathcal C$ is the speed at which its risk converges to zero when the number $n$ of available observations tends to infinity. For all the estimators defined in the sequel, we are interested in upper bounds on their risk, in order to get information about their speed of convergence. The minimax risk on a class $\mathcal X$, when $n$ observations are available, is defined as
\begin{equation*}
	\mathcal R_n(\mathcal C)=\inf_{\hat G_n} \mathcal R_n(\hat G_n ; \mathcal C), \tag{$**$}
\end{equation*}
where the infimum is taken over all set estimators depending on $n$ observations. If $\mathcal R_n(\mathcal C)$ converges to zero, we call the minimax rate of convergence on the class $\mathcal C$ the speed at which $\mathcal R_n(\mathcal C)$ tends to zero.
It is interesting to provide a lower bound for $\mathcal R_n(\mathcal C)$. By definition, no estimator can achieve a better rate on $\mathcal C$ than that of the lower bound. This bound gives also information on how close the risk of a given estimator is to the minimax risk. If the rate of the upper bound on the risk of an estimator matches the rate of the lower bound on the minimax risk on the class $\mathcal C$, then the estimator is said to have the minimax rate of convergence, or to be optimal in the minimax sense - up to constants - on this class.

A test consists of deciding whether to reject or not a given hypothesis, called the null hypothesis, when it is compared to an alternative one. Let $h\in (0,1)$. In the whole paper, we will consider the following null hypothesis:
	$$H_0 : G=\emptyset,$$
and the alternative hypothesis:
	$$H_1 : |G|\geq h.$$
Testing $H_0$ against $H_1$ is equivalent to deciding whether there is a set $G$ or not. A test is a random variable $\tau_n$, which is built from the data set, and whose possible values are 0 and 1. The decision associated to the test $\tau_n$ is to reject $H_0$ if and only if $\tau_n=1$. We measure the performance of a test $\tau_n$ on a class $\mathcal C$ using 
$$\gamma_n(\tau_n,\mathcal C)=\PP_\emptyset\left[\tau_n=1\right]+\sup_{G\in\mathcal C,|G|\geq h}\PP_G\left[\tau_n=0\right].$$
This quantity is the sum of the errors of the first and the second kinds of the test $\tau_n$.
We say that $\tau_n$ is consistent on the class $\mathcal C$ if and only if $\gamma_n(\tau_n,\mathcal C)\longrightarrow 0$, when $n\rightarrow\infty$. Let us allow the number $h$ to depend on $n$. We call the separation rate on the class $\mathcal C$ any sequence of positive numbers $r_n$ such that:
\begin{itemize}
	\item if $\displaystyle{\frac{h}{r_n} \underset{n\rightarrow\infty}{\longrightarrow}\infty}$, then there exists a consistent test on $\mathcal C$, and
	\item if $\displaystyle{\frac{h}{r_n} \underset{n\rightarrow\infty}{\longrightarrow} 0}$, then no test is consistent test on $\mathcal C$.
\end{itemize}

With regard to both Assumptions \ref{A1} and \ref{A1}, we focus on three different classes of sets, which are defined as bellow: 
\begin{description}
	\item[- ] $\mathcal S=\left\{[a,b] : 0\leq a\leq b\leq 1\right\}$ is the class of all segments on $[0,1]$,
	\item[- ] $\mathcal S_0=\left\{[0,\theta] : 0\leq \theta\leq 1\right\}$ is the class of all segments on $[0,1]$, satisfying Assumption \ref{A1},
	\item[- ] $\mathcal S(\mu)=\left\{G\in\mathcal S : |G|\geq \mu\right\}$ is the class of all segments on $[0,1]$, satisfying Assumption \ref{A2}.
\end{description}

In the whole paper, we assume that the noise terms $\xi_1,\ldots,\xi_n$ from Model \eqref{Model} satisfy 
\begin{equation}\label{subgauss}
		\E\left[e^{u\xi_i}\right]\leq e^{\sigma^2 u^2/2}, \forall u\in\R,
\end{equation}
	for some positive constant $\sigma>0$. This constant need not be known. By satisfying this inequality, the noise is said to be subgaussian. An important case of such subgaussian variables is zero mean Gaussian variables, for which $\sigma^2$ turns to be the variance.
	
  Since the design and the noise are assumed to be independent, reordering the $X_i$'s does not modify the model. Indeed, there exists a reordering $\{i_1,\ldots,i_n\}$ of $\{1,\ldots,n\}$, such that $X_{i_1}\leq\ldots\leq X_{i_n}$. The random indexes $i_1,\ldots,i_n$ are independent of the noise, and therefore the new noise vector $(\xi_{i_1}\leq\ldots\leq \xi_{i_n})$ has the same distribution as $(\xi_1\leq\ldots\leq \xi_n)$. With regard to this remark, we assume from now on that $\mathcal X$ is the reordering of a preliminary design, and therefore $X_1\leq\ldots\leq X_n$ almost surely, without loss of generality.

For two real valued sequences $A_n$ and $B_n$, and a parameter $\vartheta$, which may be multidimensional, we will write $A_n\asymp_\vartheta B_n$ to say that there exist positive constants $c(\vartheta)$ and $C(\vartheta)$ which depend on $\vartheta$, such that $c(\vartheta)B_n\leq A_n\leq C(\vartheta)B_n$, for $n$ large enough. If we put no subscript under the sign $\asymp$, this means that the involved constants are universal, i.e. do not depend on any parameter.

\subsection{Contributions of the paper}
	In this paper, we propose three new results. 
	\begin{enumerate}
		\item The first result is a preliminary to the rest of the paper, and has been already mostly covered in the existing literature. Model \eqref{Model} is a simple change point problem. If Assumption \ref{A1} is satisfied, the change point is unique. Otherwise, it is known that there are two change points, and the size of the jump is also known. In a hypothesis testing setup, we give the the minimal size - as a function of $n$ - of the unknown set $G$, in Model \eqref{Model}, must satisfy, for that set to be detectable, depending on Assumption \ref{A1} being satisfied or not. We prove our results under both designs (DD) and (RD).
		\item The second result concerns the change point problem is not necessary. We show that the speed $1/n$ of convergence in estimating the breakpoint $\theta$ in \eqref{ChangePoint} can be achieved without assuming $a\leq \theta\leq b$ for some positive numbers $0<a\leq b<1$.
		\item The last result shows that recovering a set requires that set to be sufficiently big, in order to achieve the parametric rate $1/n$ uniformly. We show that without this assumption, a logarithmic factor appears in the minimax rate.
	\end{enumerate}
Our work is made in the minimax setup. We provide exact rates of convergence of the minimax risks on the classes which we study. These rates are obtained by precise deviation inequalities for all the estimators that we propose. In the second and third results, the minimax rates correspond exactly to the separation rates given by the first result, for testing emptiness of the set $G$.
	
\subsection{Organization of the paper}
	This paper is organized as follows. Next section deals with the question of deciding whether there is a set or not. We give the separation rates for the corresponding detection test, for the change point problem - i.e. under Assumption \ref{A1} - and on the class $\mathcal S$. The third section deals with estimating the unknown set, and we give the minimax rates of estimation on the classes $\mathcal S$, $\mathcal S_0$ and $\mathcal S(\mu)$. The forth section is dedicated to the conclusion, followed by a discussion about possible extensions. The fifth and last section is devoted to the proofs.

\section{Detection of a set}
	We test the null hypothesis $H_0 : G=\emptyset$, against the alternative hypothesis $H_1 : |G|\geq h$. We find the asymptotic minimal value of $h$, as a function of $n$, so the two hypotheses are well separated and there exists a consistent test. This minimal value depends on the class of possible sets $G$. Intuitively, the \textit{a priori} knowledge that the unknown set $G$ belongs to the class $\mathcal S_0$ gives an important piece of information about the location of this set, and therefore it should be easier to detect it. Actually, the following theorem confirms this intuition, by showing that the separation rate is smaller - by a logarithmic factor - for the subclass $\mathcal S_0$ than for the whole class $\mathcal S$.

The idea, for the class $\mathcal S_0$, is the following. Under $H_1$, $[0,h]\subseteq G$. Therefore, we check among those observations $(X_i,Y_i)$ for which $X_i\leq h$ if there is a sufficiently large number of $Y_i$'s that are large, e.g. larger than $1/2$. Let $N=\max\{i=1,\ldots,n : X_i\leq h\}=\#\left(\mathcal X\cap [0,h]\right)$. Let $S$ be the following test statistics:
	\begin{equation*}
		S=\#\{i=1,\ldots,N : Y_i\leq \frac{1}{2}\}.
	\end{equation*}
If the alternative hypothesis holds, i.e. if $|G|\geq h$, all the $X_i, i=1,\ldots,N$ fall inside the set $G$, and the corresponding $Y_i$ should not be too small. The test statistic $S$ counts how many of these $Y_i$'s are suspiciously small. This is how is built the test $T_n^0$: 
$$T_n^0=\mathds 1(S\leq cN),$$
where $c\in\left(\PP[\xi_1\leq -1/2]),\PP[\xi_1\leq 1/2]\right).$

For the class $\mathcal S$, we propose a scan test, i.e., a procedure which scans the whole frame $[0,1]$ and seeks for a large enough quantity of successive observations for which $Y_i$ is large. If $G\in\mathcal S$, let $R(G)=\sum_{i=1}^n Y_i\mathds 1(X_i\in G)-\frac{\#(\mathcal X\cap G)}{2}$, and $R=\sup_{|G|\geq h}R(G)$. Under the alternative hypothesis, $R$ should be quite large, and we define the test $T_n^1=\mathds 1(R\geq 0)$.

	\begin{theorem} \label{TestThm}
		Let Model \eqref{Model} hold. 
		\begin{enumerate} 
			\item Assume that the design is (DD) or (RD), and that the noise satisfies:
						$$\PP[\xi_1\leq -1/2] < \PP[\xi_1\leq 1/2].$$ 
						Then, if $nh\longrightarrow\infty$, the test $T_n^0$ is consistent, i.e. $\gamma_n(T_n^0,\mathcal S_0)\longrightarrow 0$. If, in addition, the noise is Gaussian, then a separation rate on the class $\mathcal S_0$ is $r_n=1/n$.
			\item Assume that the design is (DD) or (RD). Then, if $nh/\ln n\longrightarrow \infty$, the test $T_n^1$ is consistent, i.e. $\gamma_n(T_n^1,\mathcal S)\longrightarrow 0$. If, in addition, the noise is Gaussian, then a separation rate on the class $\mathcal S$ is $r_n=\ln(n)/n$.
		\end{enumerate}
	\end{theorem}

In the next section, we show that the separation rates given in Theorem \ref{TestThm} are the minimax rates of convergence on the corresponding classes.

\section{Estimation of a set}

\subsection{Least square estimators}\label{Least square estimators}

	Let Model \eqref{Model} hold. For $G'\in\mathcal S$, let $A_0(G')=\sum_{i=1}^n\left(Y_i-\mathds 1(X_i\in G')\right)^2$ be the sum of squared errors. A way to estimate $G$ is to find a random set $\hat G_n$ which minimizes $A_0(G')$, among all possible candidates $G'$. Note that minimizing $A_0(G')$ is equivalent to maximizing 
\begin{equation}
	\label{LSECriterion}
	A(G')=\sum_{i=1}^n(2Y_i-1)\mathds 1(X_i\in G').
\end{equation}
Denote by $S=\{i=1,\ldots,n : X_i\in G\}$ and by $S'=\{i=1,\ldots,n : X_i\in G'\}$, for some $G'\in\mathcal S$. Denote by $\#(\cdot)$ the cardinality, for finite sets. The criterion $A(G')$ becomes, if denoted as a function of $S'$,
	\begin{align*}
		A(S') & = \sum_{i\in S'}(2Y_i-1) \\
		& = \sum_{i\in S'}\left(2\mathds 1(X_i\in G)+2\xi_i-1\right) \\
		& = 2\#(S\cap S')-\#S'+2\sum_{i\in S'}\xi_i,
	\end{align*}
so, 
	\begin{equation}
		\label{GausProc}
		A(S')-A(S) = -\#(S\triangle S')+2\left(\sum_{i\in S'\backslash S}\xi_i-\sum_{i\in S\backslash S'}\xi_i\right).
	\end{equation}
A subset $S'$ of $\{1,\ldots,n\}$ is called convex if and only if it is of the form $\{i,\ldots,j\}$, for some $1\leq i\leq j\leq n$. It is clear that if a convex subset $S'$ of $\{1,\ldots,n\}$ maximizes $A(S')-A(S)$, then by defining $G'=[X_{\min S'},X_{\max S'}]$, the segment $G'$ maximizes $A(G')$ (cf. \eqref{LSECriterion}).

\subsection{Estimation of the change point}
	Let Model \eqref{Model} hold, with design (DD). Assume that $G$ belongs to $\mathcal S_0$. This is the change point problem. For some $\theta\in[0,1]$, $G$ can be written as $G=[0,\theta]$. Let us make one preliminary remark. For any estimator $\hat G_n$ of $G$, the random segment $\tilde G_n= [0,\sup \hat G_n]$ performs better than $\hat G_n$, since $|\tilde G_n\triangle G|\leq|\hat G_n\triangle G|$ almost surely. Therefore, it is sufficient to consider only estimators of the form $\hat G_n=[0,\hat\theta_n]$, where $\hat \theta_n$ is a random variable. Then, $|\hat G_n\triangle G|=|\hat\theta_n-\theta|$, and the performance of the estimator $\hat G_n$ of $G$ is that of the estimator $\hat \theta_n$ of $\theta$.
Let us build a least square estimator (LSE) of $\theta$. For $M=1,\ldots,n$, let 
	\begin{align*}
		F(M) & =A(\{1,\ldots,M\} \\
		& = \sum_{i=1}^M(2Y_i-1).
	\end{align*}
Let $\hat M_n\in\underset{M=1,\ldots,n}{\operatorname{ArgMax}} \text{ } F(M)$, and $\hat\theta_n=X_{\hat M_n}$. The following theorem follows.

	\begin{theorem}\label{ThmChangePointUB}
		Let $n\geq 1$. Let Model \eqref{Model} hold, with design (DD). Let $\hat G_n=[0,\hat \theta_n]$. Then, 
		\begin{equation*}
			\sup_{G\in\mathcal S_0}\PP_G\left[|\hat G_n\triangle G|\geq \frac{x}{n}\right]\leq C_0e^{-x/(8\sigma^2)}, \forall x>0,
		\end{equation*}
		where $C_0$ is a positive constant which depends on $\sigma$ only.
	\end{theorem}

A simple application of Fubini's theorem leads to the following result.

\begin{corollary}
	Let the assumptions of Theorem \ref{ThmChangePointUB} be satisfied. Then, for all $q>0$, there exists a positive constant $A_q$ which depends on $q$ and $\sigma$ only, such that
	\begin{equation*}
		\sup_{G\in\mathcal S_0} \E_G\left[|\hat G_n\triangle G|^q\right] \leq \frac{A_q}{n^q}.
	\end{equation*}
\end{corollary}

This corollary shows that the minimax risk on the class $\mathcal S_0$ is bounded from above by $1/n$, up to multiplicative constants. Next theorem proves that up to constants, $1/n$ is also a lower bound on the minimax risk, if the noise is Gaussian. 

\begin{theorem}\label{ThmChangePointLB}
	Consider Model \eqref{Model}, with design (DD). Then for all integer $n\geq 1$,
	\begin{equation*}
		\mathcal R_n(\mathcal S_0)\geq \frac{1}{2n}.
	\end{equation*}
\end{theorem}

Combining Theorems \ref{ThmChangePointUB} and \ref{ThmChangePointLB} yields 

\begin{theorem}
	Consider Model \eqref{Model}, with design (DD). Then, the minimax risk on the class $\mathcal S$ satisfies 
	\begin{equation*}
		\mathcal R_n(\mathcal S_0)\asymp_\sigma \frac{1}{n}.
	\end{equation*}
\end{theorem}

\subsection{Recovering any set}
	Let us now assume that the unknown set $G$ does not necessarily contain 0. We shall prove that whether to assume that $G$ belongs to the class $\mathcal S(\mu)$ or not does not lead to the same minimax rate. 
	As we saw in Section \ref{Least square estimators}, an estimator of $G$ in Model \eqref{Model} can be obtained by maximizing the Gaussian process \eqref{GausProc} over all segments of $\{1,\ldots,n\}$. This is not the track that we will borrow, but it would be interesting to work precisely on this process. This would probably be the first step to extensions of our results in higher dimensions. However, this problem remains open for now. The methods that we develop in this section are quite different. If $G$ is only assumed to belong to the biggest class $\mathcal S$, the proposed estimator is the LSE, which was already detailed in \cite{Brunel2013} for convex polytopes, in higher dimension. If $|G|$ is \textit{a priori} known to be greater or equal to $\mu$, then we first build a preliminary estimator of $G$ - the LSE -, using one half of the observed sample. This estimator is not optimal, but it is close to $G$ with high probability. We show that the middle point $\hat m_n$ of this estimator is in $G$ with high probability. This brings us back to the change point problem, where $0$ is now replaced by $\hat m_n$, and we use the second half of the observed sample to estimate two change points.
	
	Let us first state the following theorem, which is, for the design (RD), a particular case of \cite[Theorem 1]{Brunel2013}, for $d=1$. 
	\begin{theorem}\label{LSETheorem1}
		Let $n\geq 2$. Let Model {Model} hold, with design (DD) or (RD).
		Let $\hat G_n\in\underset{G'\in\mathcal S}{\operatorname{ArgMax}} \text{ } A(G')$ be a LSE estimator of $G$. Then, there exist two positive constants $C_1$ and $C_2$ which depend on $\sigma$ only, such that
		\begin{equation*}
			\sup_{G\in\mathcal S} \PP_G\left[n\left(|\hat G_n\triangle G|-\frac{4\ln n}{C_2 n}\right)\geq x\right]\leq C_1e^{-C_2 x}, \forall x>0.
		\end{equation*}
	\end{theorem}
	
The expressions of $C_1$ and $C_2$ are given in the proof of \cite[Theorem 1]{Brunel2013}, for the design (RD). For the design (DD), we do not give a proof of this theorem here, but it can be easily adapted from that of the case of the design (RD).
The next corollary is immediate. 

\begin{corollary}
	Let the assumptions of Theorem \ref{LSETheorem1} be satisfied. Then, for all $q>0$, there exists a positive constant $B_q$ which depends on $q$ and $\sigma$ only, such that
	\begin{equation*}
		\sup_{G\in\mathcal S} \E_G\left[|\hat G_n\triangle G|^q\right] \leq B_q\left(\frac{\ln n}{n}\right)^q.
	\end{equation*}
\end{corollary}

This corollary shows that the minimax risk on the class $\mathcal S$ is bounded from above by $\ln(n)/n$, up to a multiplicative constant. The following theorem establishes a lower bound, if the noise is supposed to be Gaussian.

\begin{theorem}\label{thmLBBru13}
	Consider Model \eqref{Model}, with design (DD) or (RD). Assume that the noise terms $\xi_i$ are i.i.d. Gaussian random variables, with variance $\sigma^2>0$. For any large enough $n$, 
	\begin{equation*}
		\mathcal R_n(\mathcal S)\geq \frac{\alpha^2\sigma^2\ln n}{n},
	\end{equation*}
	where $\alpha$ is a universal positive constant.
\end{theorem}

This lower bound comes from \cite[Theorem 2]{Brunel2013} in the case of the design (RD), and the proof is easily adapted for the design (DD). Eventually, the minimax risk on the class $\mathcal S$ is of the order $\ln(n)/n$:

\begin{theorem}
		Consider Model \eqref{Model}, with design (DD) or (RD). Assume that the noise terms $\xi_i$ are i.i.d. Gaussian random variables, with variance $\sigma^2>0$. The minimax risk on the class $\mathcal S$ satisfies, asymptotically:
		$$\mathcal R_n(\mathcal S) \asymp_\sigma \frac{\ln n}{n}.$$
\end{theorem}

For the design (DD), we combine both Theorems \ref{ThmChangePointUB} and \ref{LSETheorem1} to find the minimax rate on the class $\mathcal S(\mu)$. 
Let Model \eqref{Model} hold, and let $G\in\mathcal S(\mu)$. First, we split the sample into two equal parts. 
Let $\mathcal D_0$ be the set of sample points with even indexes, and $\mathcal D_1$ the set of sample points with odd indexes. Note that $\mathcal D_0\cup\mathcal D_1$ is exactly the initial sample, that these two subsample are independent, and that each of them is made of at least $(n-1)/2$ data. Let $\hat G_n$ be the LSE estimator of $G$ given in Theorem \ref{LSETheorem1}, built from the subsample $\mathcal D_0$. Let $\hat m_n$ be the middle of $\hat G_n$. As it will be shown in the proof of the next theorem, $\hat m_n$ satisfies both following properties, with high probability:
\begin{enumerate}
	\item $\hat m_n\in G$,
	\item $\mu/2\leq\hat m_n\leq 1-\mu/2$.
\end{enumerate}
This brings us to estimating two change points - the endpoints of $G$ - , using the second subsample $\mathcal D_1$. 
From Theorem \ref{ThmChangePointUB}, we know that this can be done at the speed $1/n$, up to multiplicative constants. 

\begin{theorem}\label{theoremMu}
		Consider Model \eqref{Model}, with design (DD). There exists an estimator $\tilde G_n$ of $G$, such that
				\begin{equation*}
					\sup_{G\in\mathcal S(\mu)}\PP_G\left[|\tilde G_n\triangle G|\geq \frac{x}{n}\right] \leq 2C_0e^{-\mu x/(256\sigma^2)}+C_1n^4e^{-C_2\mu n/2}, \forall x>0,
				\end{equation*}
		for $n$ large enough. The positive constants $C_0$ and $C_2$ appeared in Theorems \ref{ThmChangePointUB} and \ref{LSETheorem1}.
\end{theorem}

Naturally, Theorem \ref{theoremMu} leads to the next corollary.

\begin{corollary}
	Let the assumptions of Theorem \ref{theoremMu} be satisfied. Then, for all $q>0$, there exists a positive constant $B'_q$ which depends on $q$, $\mu$ and $\sigma$ only, such that
	\begin{equation*}
		\sup_{G\in\mathcal S(\mu)}\E_G\left[|\tilde G_n\triangle G|^q\right] \leq \frac{B'_q}{n^q}.
	\end{equation*}
\end{corollary}

This corollary, for $q=1$, shows that the minimax risk on the class $\mathcal S(\mu)$ is bounded from above by $1/n$, up to a multiplicative constant. A very similar proof to that of Theorem \ref{ThmChangePointLB} yields a lower bound for this minimax risk, which leads to the next theorem.

\begin{theorem}
		Consider Model \eqref{Model}, with design (DD). The minimax risk on the class $\mathcal S(\mu)$ satisfies, asymptotically,
		\begin{equation*}
			\mathcal R_n(\mathcal S(\mu))\asymp_{\mu,\sigma} \frac{1}{n}.
		\end{equation*}
\end{theorem}

\begin{remark}
	Note that, in Theorem \ref{theoremMu}, the upper bound contains one residual term which does not depend on $x$. This term, in order to be sufficiently small, requires that $\mu$ - if allowed to depend on $n$ - is of larger order than $\ln(n)/n$. This in an echo to Theorem \ref{TestThm}, in which we showed that the smallest set which can be detected has measure of this order exactly. In addition, if $\mu$ is of the order of $\ln(n)/n$, then the proof of the lower bound of Theorem \ref{thmLBBru13} can be applied, and the minimax risk on the class $\mathcal S(\mu)$ will satisfy $\mathcal R_n(\mathcal S(\mu))\asymp_{\sigma} \frac{\ln n}{n}$.
\end{remark}

\section{Conclusion and discussion}

We summarize our results in Table \ref{Table1}. The rates that are written in this table hold for Gaussian noise, which is the most important case. In each case we indicate the design for which the rate holds. \\

\begin{table}[h]
\centering
\begin{tabular}{|c|c|c|c|}
\hline
& $\mathcal S$ & $\mathcal S_0$ & $\mathcal S(\mu)$ \\
\hline
Minimax rate & $\ln(n)/n$ (DD,RD) & $1/n$ (DD) & $1/n$ (DD) \\
\hline
Separation rate & $\ln(n)/n$ (DD,RD) & $1/n$ (DD,RD) & $\cdot$ \\
\hline
\end{tabular}
\caption{Minimax risks and separation rates for the classes $\mathcal S$, $\mathcal S_0$ and $\mathcal S(\mu)$.}
\label{Table1}
\end{table}

Note that in two cases, only the design (DD) has been considered. This is mainly for technical reasons, and we believe that the rates are still the same for the design (RD).

It comes out that asymptotically, a segment can be estimated infinitely faster when it is \textit{a priori} supposed either to contain a given point (here, $0$), or to be large enough. The main question that remains is: does this phenomenon still hold for two - or higher - dimensional sets ? An important - if not essential - assumption which has been done all over this paper is the convexity of the unknown set. In dimension 1, the class of convex subsets of $[0,1]$ is simple, and parametric. In any higher dimension $d$, the class $\mathcal C_d$ of convex subsets of $[0,1]^d$ is much more complex. In particular, its metric entropy is much larger than that of a parametric family \cite{Bronshtein}, and it seems to us that it is the complexity of this class that makes it harder to estimate a set, than detectability. 

In Model \eqref{Model}, if $G$ belongs to $\mathcal C_d$, estimation of $G$ can be done at the minimax rate $n^{-2/(d+1)}$ \cite{Brunel2013}. However, an adaptation of the proof of \cite[Theorem 2]{Brunel2013} shows the following: on any subclass of $\mathcal C_d$ invariant under translations and invertible affine transformations - which keep a set $G$ inside the frame $[0,1]^d$ -, the minimax rate is at least of order $\ln(n)/n$. In addition, we believe this is the separation rate for the detection problem, on any such subclass of $\mathcal C_d$. 

For the whole class $\mathcal C_d$, since $\ln(n)/n$ is much smaller than $n^{-2/(d+1)}$, the minimax rate is of the order of $n^{-2/(d+1)}$. However, for a parametric subclass, such as that of all convex polytopes with a given number of vertices, we believe that the complexity of the class leads to a term of order $1/n$ in the minimax risk, although detectability leads to a term of order $\ln(n)/n$. In our opinion, this is what explains that, as shown in \cite{Brunel2013}, the minimax risk on the class of all convex polytopes of $[0,1]^d$, with a given number of vertices, is of order $\ln(n)/n$. 

Motivated by Theorem \ref{theoremMu}, we also conjecture that $1/n$ is the minimax rate on the class of all convex polytopes of $[0,1]^d$, of volume greater than a given $\mu>0$. However, it is not possible to extend Theorem \ref{theoremMu} to higher dimensions. An adaptation of the proof of \cite[Theorem 2]{Brunel2013}, by taking, as the $M$ hypotheses used in the proof, sets which contain the origin and have pairwise zero measure intersections, would show that $(\ln n)/n$ remains a lower bound under that assumption. Yet, we believe that if the unknown set is assumed to contain a given section of positive $d-1$ dimensional Lebesgue measure, of a given hyperplane in $\R^d$, then an analog of Theorem \ref{theoremMu} should hold, and the minimax rate should be of order $1/n$.

To end this discussion, we note that Model \eqref{Model2} can be rewritten in the vector form:
\begin{equation}
	y=\beta+\xi,
\label{vectorform0}
\end{equation}
where $y=(Y_1,\ldots,Y_n)^\top$, $\beta=\left(\delta\mathds 1(X_1\in G),\ldots,\delta\mathds 1(X_n\in G)\right)^\top \in \{0,\delta\}^n$ and $\xi=(\xi_1,\ldots,\xi_n)^\top$. 
If we denote by $A$ the $n\times n$ lower triangular matrix with coefficients on and under the main diagonal equal to $1$, equation \eqref{vectorform0} can be rewritten as:
\begin{equation}
	y=A\gamma+\xi,
\label{vectorform1}
\end{equation}
where $\gamma=(\beta_1,\beta_2-\beta_1,\ldots,\beta_n-\beta_{n-1})^\top$ is the vector of differences of $\beta$. The vector $\gamma$ is sparse: it has exactly two nonzero coefficients, corresponding to the change points locations. This formulation of our initial problem leads to a high dimensional regression setup, under sparsity constraints. If $\delta$ is known, equal to $1$ -this is Model \eqref{Model}-, we are not interested in the estimation of the vector $\gamma$ itself, but only of its support, which indicates the location of $G$, if the design is (DD). If $\delta$ is unknown, the two nonzero coefficients of $\gamma$ take values $\delta$ and $-\delta$, and it is also of interest to estimate $\gamma$ itself, in addition to its support. Penalized regression methods, with Lasso-type penalizations (see \cite{Tibshirani1996} for details about Lasso), have been addressed in \cite{MammenVanDeGeer1997}, \cite{BoysenKempeLiebscherMunkWittich2009}, among others. In these two works, the penalization is written in terms of the number of jumps in the vector $\beta$ in \eqref{vectorform0}, which is equal, for Model \eqref{Model}, to the number of nonzero coefficients of $\gamma$ in \eqref{vectorform1}. Therefore, the penalty is equivalent to a $L_0$ or BIC-type penalty, see for example \cite{BuneaTsybakovWegkamp2007}. Estimation of the support is usually a secondary question, which is addressed in terms of consistency, and not of rate of convergence. The main focus is about estimation of $\gamma$, or prediction, i.e., estimation of $\beta$. We believe that the parameter $\mu$, if both $\mu$ and $G$ are unknown, can be estimated at the speed $\sqrt{(\ln n)/n}$. In that case, some results about estimation rates of $G$, for the design (DD), are given in \cite{BoysenKempeLiebscherMunkWittich2009}, but the risks, in expectation, are not computed. 

It would be interesting to understand how Lasso or BIC-type estimators could be adapted in order to achieve the rate $1/n$ on the classes $\mathcal S_0$ and $\mathcal S(\mu)$. The penalties of the classical estimators do not use the information that is contained in those classes, and we believe that these estimators are suboptimal. However, one should try to modify the penalties, according to that information, in order to improve the accuracy of these estimators.

\section{Proofs}

\subsection{Proof of Theorem \ref{TestThm}}

\paragraph{On the class $\mathcal S_0$}

\subparagraph{Upper bound}
	Let us first prove the upper bound, i.e. assume that $nh\rightarrow\infty$, and prove that there exists a consistent test. Recall that $N=\max\{i=1,\ldots,n : X_i\leq h\}=\#\left(\mathcal X\cap [0,h]\right)$. If the design is (DD), then $N$ is just equal to the integer part of $nh$. If the design (RD), then $N$ is a binomial random variable, with parameters $n$ and $h$.
Let us show first that the error of the first kind of the test $T_n^0$ goes to zero, when $n\rightarrow\infty$.
\begin{align*}
	\PP_\emptyset\left[S\leq cN\right] & = \PP_\emptyset\left[\#\left\{i=1,\ldots,N :Y_i>\frac{1}{2}\right\}\geq(1-c)N\right] \\
	& \leq \E\left[\PP_\emptyset\left[\#\left\{i=1,\ldots,N : \xi_i>\frac{1}{2}\right\}\geq(1-c)N|\mathcal X\right]\right].
\end{align*}
Since the $\xi_i$'s are independent of $\mathcal X$, the distribution of $\#\{i=1,\ldots,N : \xi_i>\frac{1}{2}\}$ conditionally to $\mathcal X$ is binomial, with parameters $N$ and $\beta$, where $\beta =\PP\left[\xi_1>1/2\right]\in [0,1)$.
Thus, by Bernstein's inequality for binomial random variables, by defining $\displaystyle{\gamma=\frac{(1-c-\beta)^2}{2\beta(1-\beta)+(1-c-\beta)/3}>0}$,
\begin{equation*}
	\PP_\emptyset\left[S\leq cN\right] \leq \E\left[\exp\left(-\gamma N\right)\right].
\end{equation*}
If $\mathcal X$ satisfies (DD), then $N\geq nh-1$ and it is clear that $\PP_\emptyset\left[S\leq cN\right]\longrightarrow 0$.
If $\mathcal X$ satisfies (RD), then 
\begin{equation*}
	\E\left[\exp\left(-\gamma N\right)\right]=\exp\left(-nh\left(1-e^{-\gamma}\right)\right),
\end{equation*} 
so $\PP_\emptyset\left[S\leq cN\right]\longrightarrow 0$.

Let us show, now, that the error of the second kind goes to zero as well. Let $G\in\mathcal S_0$ satisfying the alternative hypothesis, i.e. $|G|\geq h$. Denote by $\beta'=\PP[\xi_1\leq -1/2]$, and by $\displaystyle{\gamma'=\frac{(c-\beta')^2}{2\beta'(1-\beta')+(c-\beta')/3}>0}$
\begin{align*}
	\PP_G\left[S>cN\right] & = \PP_G\left[\#\left\{i=1,\ldots,N : Y_i\leq\frac{1}{2}\right\}> c N\right] \\
	& \leq \E\left[\PP_\emptyset\left[\#\left\{i=1,\ldots,N : \xi_i\leq-\frac{1}{2}\right\}>cN|\mathcal X\right]\right] \\
	& \leq \E\left[\exp\left(-\gamma' N\right)\right],
\end{align*}
by a similar computation to that for the error of the first kind. Since the right-side of the last inequality does not depend on $G$, 
\begin{equation*}
	\sup_{|G|\geq h}\PP_G\left[S>cN\right] \leq \E\left[\exp\left(-\gamma' N\right)\right]
\end{equation*}
and therefore, by the same argument as for the error of the first kind, goes to zero when $n\rightarrow\infty$, for both designs (DD) and (RD).

\subparagraph{Lower bound}

Assume, now, that $nh\rightarrow 0$. Let $\tau_n$ be any test. Let $G_1=[0,h]$. We denote by $\mathcal H$ the Hellinger distance between probability measures. The following computation uses properties of this distance, which can be found in \cite{Tsybakov2009}. 
\begin{align}
	\gamma_n(\tau_n,\mathcal C) & \geq \E_\emptyset\left[\tau_n\right]+\E_{G_1}\left[1-\tau_n\right] \nonumber \\
	& \geq \int \min\left(d\PP_\emptyset,d\PP_{G_1}\right) \nonumber \\
	& \geq \frac{1}{2}\left(1-\frac{\mathcal H(\PP_\emptyset,\PP_{G_1})}{2}\right)^2. \label{LBTest1}
\end{align}
Let $G,G'\in\mathcal S$. A simple computation shows that, for the design (DD),
\begin{equation}\label{hellinger1}
	1-\frac{\mathcal H(\PP_G,\PP_{G'})}{2} = \exp\left(-\frac{\#\left(\mathcal X\cap (G\triangle G')\right)}{8\sigma^2}\right),
\end{equation}
and for the design (RD), 
\begin{equation}\label{hellinger2}
	1-\frac{\mathcal H(\PP_G,\PP_{G'})}{2} = \left(1-\left(1-e^{-\frac{1}{8\sigma^2}}\right)|G\triangle G'|\right)^n.
\end{equation}
In particular, for the design (DD), 
\begin{equation*}
	1-\frac{\mathcal H(\PP_\emptyset,\PP_{G_1})}{2} \geq \exp\left(-\frac{nh}{8\sigma^2}\right),
\end{equation*}
and for the design (RD),
\begin{equation*}
	1-\frac{\mathcal H(\PP_\emptyset,\PP_{G_1}^{\otimes n})}{2} = \left(1-\left(1-e^{-\frac{1}{8\sigma^2}}\right)h\right)^n.
\end{equation*}
In both cases, we showed that the right side of \eqref{LBTest1} tends to $1/2$, when $n\rightarrow\infty$. Therefore the test $\tau_n$ is not consistent.

\paragraph{On the class $\mathcal S$}

\subparagraph{Upper bound}

Assume that $\displaystyle{\frac{nh}{\ln n}\longrightarrow\infty}$. Let us first show that the error of the first kind of $T_n^1$ goes to zero, when $n\rightarrow\infty$. Recall that $T_n^1=\mathds 1(R\geq 0)$, where $R=\sup_{|G|\geq h}R(G)$ and $R(G)=\sum_{i=1}^n Y_i\mathds 1(X_i\in G)-\frac{\#(\mathcal X\cap G)}{2}$, for all $G\in\mathcal S$. Note that $R(G)$ is piecewise constant, and can only take a finite number of values. It is clear that 
\begin{equation*}
	\{R(G):G\in\mathcal S,|G|\geq h\}=\{R([X_k,X_l)):1\leq k< l\leq n, X_l-X_k>h\}.
\end{equation*}
Recall that for $1\leq k< l\leq n$, $\displaystyle{R([X_k,X_l))=\frac{1}{2}\sum_{i=k}^{l-1}(2Y_i-1)}$.
Therefore, 
\begin{align}
	\PP_\emptyset[R\geq 0] & = \PP_\emptyset\left[\max_{\substack{1\leq k < l\leq n \\ X_l-X_k>h}} R([X_k,X_l))>0\right] \nonumber \\
	& \leq \PP_\emptyset\left[\bigcup_{1\leq k < l\leq n} \left\{R([X_k,X_l))>0\right\}\cap\left\{X_l-X_k>h\right\}\right] \nonumber \\
	& \leq \sum_{1\leq k < l\leq n}\PP_\emptyset\left[R([X_k,X_l))>0,X_l-X_k>h\right] \nonumber \\
	& \leq \sum_{1\leq k < l\leq n}\PP_\emptyset\left[\sum_{i=k}^{l-1}(2\xi_i-1)>0\right]\PP[X_l-X_k>h]. \label{1stErrS1}
\end{align}
For $1\leq k<l\leq n$, 
\begin{equation}
	\PP_\emptyset\left[\sum_{i=k}^{l-1}(2\xi_i-1)>0\right]\leq \exp\left(-\frac{(l-k)\sigma^2}{8}\right), \label{1stErrS2}
\end{equation}
using Markov's inequality and \eqref{subgauss}.

If the design is (DD), then $\PP[X_l-X_k>h]$ is 1 if and only if $l-k>nh$, 0 otherwise, so from \eqref{1stErrS1} and \eqref{1stErrS2} we get that
\begin{align*}
	\PP_\emptyset[R\geq 0] & \leq \sum_{l-k >nh}\exp\left(-\frac{(l-k)\sigma^2}{8}\right) \\
	& \leq \sum_{l-k >nh}\exp\left(\frac{(-nh)\sigma^2}{8}\right) \\
	& \leq \frac{n^2}{2}\exp\left(\frac{(-nh)\sigma^2}{8}\right) \longrightarrow 0,
\end{align*}
when $n\rightarrow\infty$, which proves that the error of the first kind goes to zero.

If the design is (RD), let us use the following Lemma.

\begin{lemma}\label{lemmaRD}
	Let $X_1,\ldots,X_n$ be the (RD) design. Then, for any $1\leq k<l\leq n$, and $h>0$,
	\begin{equation*}
		\PP[X_l-X_k>h] \leq  n\exp\left(-nh(1-e^{-u})+u(l-k)\right), \forall u>0.
	\end{equation*}
\end{lemma}

\paragraph{Proof of Lemma \ref{lemmaRD}}
Note that the event $\{X_l-X_k>h\}$ is equivalent to $\{\#(\mathcal X\cap (X_k,X_k+h)) < l-k\}$. Let us denote by $X_1',\ldots,X_n'$ the preliminary design, from which $X_1,\ldots,X_n$ is the reordered version. The random variables $X_1',\ldots,X_n'$ are then i.i.d., with uniform distribution on $[0,1]$. Hence,
\begin{align}
	\PP[X_l-X_k>h] & = \sum_{j=1}^n \PP\left[\#(\mathcal X\cap (X_k,X_k+h)) < l-k, X_k=X_j'\right] \nonumber \\
	& \leq \sum_{j=1}^n \PP\left[\#(\mathcal X\cap (X_j',X_j'+h)) < l-k\right] \nonumber \\
	& \leq \sum_{j=1}^n \E\left[\PP\left[\#(\mathcal X\cap (X_j',X_j'+h)) < l-k|X_j'\right]\right] \nonumber \\
	& \leq \sum_{j=1}^n \E\left[\PP\left[n-\#(\mathcal X\cap (X_j',X_j'+h)) \geq n-l+k|X_j'\right]\right] \nonumber \\
	& \leq \sum_{j=1}^n \E\left[f(X_j')\right], \label{1stErrS3}
\end{align}
where $f(x)=\PP\left[n-\#(\mathcal X\cap (x,x+h)) \geq n-l+k\right]$, for all $x\in[0,1]$.
The random variable $n-\#(\mathcal X\cap (x,x+h))$ is binomial with parameters $n$ and $1-h$, and by Markov's inequality, for all $u>0$,
\begin{equation}
	\label{1stErrS4} f(x) \leq e^{nu}\left(1-h(1-e^{-u})\right)^ne^{-u(n-l+k)},
\end{equation}
and \eqref{1stErrS3} and \eqref{1stErrS4} yield the lemma. \hfill $\Box$\\[0.4cm]

Therefore, by \eqref{1stErrS1}, \eqref{1stErrS2} and \eqref{1stErrS3}, and Lemma \ref{lemmaRD} with $u=\sigma^2/8$,
\begin{align*}
	\PP_\emptyset[R\geq 0] & \leq \sum_{1\leq k < l\leq n}\exp\left(-\frac{(l-k)\sigma^2}{8}-nh(1-e^{-u})+u(l-k)\right) \\
	& \leq \frac{n^2}{2}\exp\left(-nh(1-e^{-\sigma^2/8})\right) \longrightarrow 0,
\end{align*}
when $n\rightarrow \infty$, which proves that the error of the first kind goes to zero.

Let us bound, now, the error of the second kind. Let $G\in\mathcal S$ satisfying $|G|\geq h$. For this $G$, denote by $N_G=\#(\mathcal X\cap G)$. Then,
\begin{align}
	\PP_G[R<0] & \leq \PP_G\left[R(G)\leq N_G/2\right] \nonumber \\
	& \leq \PP\left[\sum_{i=1}^n \xi_i\mathds 1(X_i\in G)\leq -N_G/2\right]. \label{PrUB2_S}
\end{align}

For the design (DD), $N_G$ is the integer part of $n|G|$, so $N_G\geq nh$. Therefore, by Markov's inequality, and by \eqref{subgauss}, \eqref{PrUB2_S} becomes 
\begin{equation}
	\PP_G[R<0] \leq \exp\left(-\frac{N_G}{8\sigma^2}\right) \leq \exp\left(-\frac{nh}{8\sigma^2}\right). \label{DDUB2}
\end{equation}

For the design (RD), $N_G$ is a random binomial variable, with parameters $n$ and $|G|$. By conditioning to the design and using Markov's inequality, \eqref{PrUB2_S} becomes 
\begin{align}
	\PP_G[R<0] & \leq \PP\left[\sum_{i=1}^n -\xi_i\mathds 1(X_i\in G)\geq N_G/2\right] \nonumber \\
	& \leq \E\left[\exp\left(-\frac{N_G}{8\sigma^2}\right)\right] \nonumber \\
  & \leq \exp\left(-Cn|G|\right) \leq \exp\left(-Cnh\right), \label{RDUB2}
\end{align}
where $C=1-e^{-\frac{1}{8\sigma^2}}$.

In both cases \eqref{DDUB2} and \eqref{RDUB2}, the right side does not depend on $G$, and goes to zero as $n\rightarrow\infty$. We conclude that, for both designs (DD) and (RD), 
\begin{equation*}
	\sup_{|G|\geq h}\PP_G[R<0] \longrightarrow 0,
\end{equation*}
which ends the proof of the upper bound.

\subparagraph{Lower bound}

We more or less reproduce the proof of \cite{Gayraud2001}, Theorem 3.1. Here, the noise is supposed to be Gaussian, with variance $\sigma^2$. Let us assume that $\displaystyle{\frac{nh}{\ln n}\longrightarrow 0}$. Let $M=1/h$, assumed to be an integer, without loss of generality. For $q=0,\ldots,M$, let $G_q=[qh,(q+1)h]$. 
For $q=1,\ldots,M$, let $Z_q=\frac{d\PP_{G_q}}{d\PP_\emptyset}(X_1,Y_1,\ldots,X_n,Y_n)$, and denote by $\bar Z=\frac{1}{M}\sum_{q=1}^M Z_q$. Let $\tau_n$ be any test. Then,
\begin{align}
	\gamma_n(\tau_n,\mathcal S) & \geq \PP_\emptyset\left[\tau_n=1\right]+\frac{1}{M}\sum_{q=1}^M\PP_{G_q}\left[\tau_n=0\right] \nonumber \\
	& \geq \frac{1}{M}\sum_{q=1}^M\left(\PP_\emptyset\left[\tau_n=1\right]+\PP_{G_q}\left[\tau_n=0\right]\right) \nonumber \\
	& \geq \frac{1}{M}\sum_{q=1}^M\left(\E_\emptyset\left[\tau_n\right]+\E_{G_q}\left[1-\tau_n\right]\right) \nonumber \\
	& \geq \frac{1}{M}\sum_{q=1}^M\E_\emptyset\left[\tau_n+(1-\tau_n)Z_q\right] \nonumber \\
	& \geq \E_\emptyset\left[\tau_n+(1-\tau_n)\bar Z\right] \nonumber \\
	& \geq \E_\emptyset\left[\left(\tau_n+(1-\tau_n)\bar Z\right)\mathds 1(\bar Z\geq 1/2)\right] \nonumber \\
	& \geq \frac{1}{2}\PP_\emptyset\left[\bar Z\geq 1/2\right]. \label{LBtest01}
\end{align}
Let us prove that $\E_\emptyset[\bar Z]=1$, and that $\VV_\emptyset[\bar Z]\longrightarrow 0$. This will imply that the right side term of \eqref{LBtest01} goes to zero, when $n\rightarrow\infty$. 

For $q=1,\ldots,M$, under the null hypothesis, 
\begin{align}
	Z_q & = \exp\left(-\frac{1}{2\sigma^2}\sum_{i=1}^n\left((Y_i-\mathds 1(X_i\in G_q))^2-Y_i^2\right)\right) \nonumber \\
	& = \exp\left(\frac{1}{2\sigma^2}\sum_{i=1}^n(2\xi_i-1)\mathds 1(X_i\in G_q)\right). \label{Zq}
\end{align}
By its definition, $Z_q$ has expectation 1 under $\PP_\emptyset$:
\begin{equation}\label{expectZ}
	\E_\emptyset[\bar Z]=1.
\end{equation}

Since, almost surely, no design point falls in two $G_q$'s at the time, a simple computation shows that the random variables $Z_q, q=1,\ldots,M$, are not correlated. Thus,
$$\VV_\emptyset[\bar Z]=\frac{1}{M^2}\sum_{q=1}^M\VV_\emptyset[Z_q].$$
Let us bound from above $\VV_\emptyset[Z_q]$, for $q=1,\ldots,M$:
\begin{align}
	\VV_\emptyset[Z_q] & \leq \E_\emptyset[Z_q^2] \nonumber \\
	& = \E\left[\exp\left(-\frac{\#(\mathcal X\cap G_q)}{\sigma^2}\right)\E_\emptyset\left[\exp\left(\frac{2}{\sigma^2}\sum_{i=1}^n\xi\mathds 1(X_i\in G_q)\right) | \mathcal X\right]\right] \nonumber \\
	& = \E\left[\exp\left(\frac{\#(\mathcal X\cap G_q)}{\sigma^2}\right)\right]. \label{boundonvariance}
\end{align}

If the design is (DD), then we get that 
\begin{equation*}
	\VV_\emptyset[Z_q]\leq \exp\left(\frac{nh+1}{\sigma^2}\right),
\end{equation*}
and the variance of $\bar Z$ is then bounded from above:
\begin{equation}
	\VV_\emptyset[\bar Z]\leq h\exp\left(\frac{nh+1}{\sigma^2}\right). \label{var1}
\end{equation}

If the design is (RD), then $\#(\mathcal X\cap G_q)$ is a binomial random variable with parameters $n$ and $h$, so from \eqref{boundonvariance}, we get that

\begin{align*}
	\VV_\emptyset[Z_q] & \leq \left(1+\left(e^{1/\sigma^2}-1\right)h\right)^n \\
	& \leq \exp\left(Cnh\right),
\end{align*}
where $C=e^{1/\sigma^2}-1$, and the variance of $\bar Z$ is then bounded from above:
\begin{equation}
	\VV_\emptyset[\bar Z]\leq h\exp\left(Cnh\right). \label{var2}
\end{equation}

Since we assumed that $nh/\ln n\longrightarrow 0$, the right side terms of \eqref{var1} and \eqref{var2} go to zero, and therefore, for both designs (DD) and (RD),
\begin{equation}\label{boundvarianceZDD}
	\VV_\emptyset[\bar Z]\longrightarrow 0.
\end{equation} 

Finally, we get from \eqref{LBtest01}, \eqref{expectZ} and \eqref{boundvarianceZDD}, that 
\begin{equation*}
	\underset{n\rightarrow\infty}{\operatorname{liminf}}\gamma_n(\tau_n,\mathcal S)\geq \frac{1}{2}.
\end{equation*}
This concludes the proof. \hfill $\blacksquare$

\subsection{Proof of Theorem \ref{ThmChangePointUB}}

The beginning of this proof holds for any design $\{X_1,\ldots,X_n\}$, independent of the noise $\xi_i, i=1,\ldots,n$.
Let $G\in\mathcal S_0$. Let $M=\max\{i=1,\ldots,n : X_i\in G\}$ - set $M=0$ if the set is empty -. Then, $\hat M_n\in\underset{M'=1,\ldots,n}{\operatorname{ArgMax}} \text{ } \left(F(M')-F(M)\right)$, and, by \eqref{GausProc},

\begin{equation*}
	F(M')-F(M)= -|M'-M| + { \left\{
    \begin{array}{l}
        2\sum_{i=M'+1}^M\xi_i   \mbox{ }\mbox{ if } M>M',   \vspace{3mm} \\
				0 \mbox{ }\mbox{ if } M'=M, \vspace{3mm} \\
        -2\sum_{i=M+1}^{M'}\xi_i   \mbox{ }\mbox{ if } M<M'.
    \end{array}
	\right.}
\end{equation*}

Let us complete the i.i.d. sequence $\xi_1,\ldots,\xi_n$ to obtain an infinite double sided i.i.d. sequence $(\xi_i)_{i\in\mathbb Z}$, independent of the design.
Let $k\in\mathbb N^*$ be any positive integer. Define, for $i\in\mathbb Z, \tilde\xi_i=\xi_{i+M}$. Since $M$ depends on the design only, it is independent of the $\xi_i,i\in\mathbb Z$, and therefore, the $\tilde\xi_i,i\in\mathbb Z$ are i.i.d., with same distribution as $\xi_1$.
Let $E_k$ be the event $\{\hat M_n\geq M+k\}$. If $E_k$ holds, then :
\begin{equation*}
	0\leq F(\hat M_n)-F(M)=M-\hat M_n - 2\sum_{i=M+1}^{\hat M_n}\xi_i,
\end{equation*}
and it follows that 
\begin{align*}
	0 & \leq \max_{M+k\leq j\leq n}\left(M-j - 2\sum_{i=M+1}^{j}\xi_i\right) \\
	& \leq \max_{M+k\leq j\leq n}\left(M-j - 2\sum_{i=1}^{j-M}\tilde\xi_i\right) \\
	& \leq \max_{k\leq j\leq n-M}\left(-j - 2\sum_{i=1}^{j}\tilde\xi_i\right).
\end{align*}
Hence, for all $u>0$,
\begin{align*}
	\PP_G[E_k] & \leq \PP_G\left[\max_{k\leq j}\left(-j - 2\sum_{i=1}^{j}\hat\xi_i\right) \geq 0\right] \\
	& \leq \PP\left[\max_{k\leq j}\left(-j - 2\sum_{i=1}^{j}\xi_i\right) \geq 0\right] \\
	& \leq \sum_{j=k}^\infty \PP\left[-2\sum_{i=1}^{j}\xi_i \geq j\right] \\
	& \leq \sum_{j=k}^\infty \frac{\E\left[e^{-2u\sum_{i=1}^{j}\xi_i}\right]}{e^{uj}}, \mbox{ by Markov's inequality} \\
	& \leq \sum_{j=k}^\infty e^{(-u+2\sigma^2u^2)j}, \mbox{ by \eqref{subgauss}}
\end{align*}
and, by choosing $u=1/(4\sigma^2)$,
\begin{equation*}
	\PP_G[E_k]\leq Ce^{-k/(8\sigma^2)},
\end{equation*}
where $C=\left(1-e^{-1/(8\sigma^2)}\right)^{-1}$ is a positive constant.
By symmetry, we obtain that :
\begin{equation}\label{devind}
	\PP_G[|\hat M_n-M|\geq k]\leq 2Ce^{-k/(8\sigma^2)}.
\end{equation}

If the design is (DD), the conclusion is straightforward, since for all $i,j=1,\ldots,n, |X_i-X_j|=\frac{|i-j|}{n}$, and Theorem \ref{ThmChangePointUB} is proved. \hfill $\blacksquare$ \\[0.2cm]

If the design is (RD), it is not clear how to go from \eqref{devind} to an upper bound for the probability $\PP_G[|\hat \theta_n-\theta|\geq \epsilon]$, for $\epsilon>0$. This question should be addressed in a future work.

\subsection{Proof of Theorem \ref{ThmChangePointLB}}

The proof is straightforward. Let $G_1=[0,0]$ and $G_2=[0,1/(2n)]$. Then $\PP_{G_1}=\PP_{G_2}$, since no point of the design falls in $G_1\triangle G_2$, and for any estimator $\hat G_n$, 
\begin{align*}
	\sup_{G\in\mathcal S_0}\E_G\left[|\hat G_n\triangle G|\right] & \geq \E_{G_1}\left[|\hat G_n\triangle G_1|\right]+\E_{G_2}\left[|\hat G_n\triangle G_2|\right] \\
	& \geq \E_{G_1}\left[|\hat G_n\triangle G_1|+|\hat G_n\triangle G_2|\right] \\
	& \geq \E_{G_1}\left[|G_1\triangle G_2|\right] \mbox{ by the triangle inequality} \\
	& \geq \frac{1}{2n}. 
\end{align*}
\hfill $\blacksquare$

\subsection{Proof of Theorem \ref{theoremMu}}

Let $I_0$ be the set of even positive integers less or equal to $n$, and $I_1$ the set of odd such integers. Note that $\{X_i : i\in I_0\}$ is a deterministic and regular design, with step $2/n$. 
Let $\hat G_n\in\underset{G'\in\mathcal S}{\operatorname{ArgMax}} \text{ } \sum_{i\in I_0}(2Y_i-1)\mathds 1(X_i\in G')$ be the LSE estimator given in Theorem \ref{LSETheorem1}, built using only the subsample $\{X_i : i\in I_0\}$. Let $x>0$, whose value will be specified in the course of the proof. Consider the event $E_x=\{|\hat G_n\triangle G|<\frac{x\ln n}{n}\}$.
By Theorem \ref{LSETheorem1}, this event holds with probability at least $1-C_1e^{-(C_2x-4)\ln n}$.
Choose $x$ such that $\frac{x\ln n}{n}\leq \mu/2$. This choice implies that on the event $E_x$, $|\hat G_n\triangle G|<\mu\leq |G|$, so necessary, $\hat G_n$ and $G$ must intersect. Thus, still on the event $E_x$,
\begin{equation*}
	|\hat G_n\triangle G|=|\hat b_n-b|+|\hat a_n-a|,
\end{equation*}
where we denoted by $G=[a,b]$ and $\hat G_n=[\hat a_n,\hat b_n]$.
Let $m=\frac{a+b}{2}$ and $\hat m_n=\frac{\hat a_n+\hat b_n}{2}$ be, respectively, the middle points of $G$ and $\hat G_n$. From now on, let us assume that $E_x$ holds. Then,
\begin{align}
	|\hat m_n-m| & \leq \frac{1}{2}(|\hat b_n-b|+|\hat a_n-a|) \nonumber \\
	& \leq \frac{1}{2}|\hat G_n\triangle G| \nonumber \\
	& \leq \frac{x\ln n}{2n} \nonumber \\
	& \leq \frac{\mu}{4}. \label{middles}
\end{align}
Therefore, $\hat m_n\in G$ and, combining \eqref{middles} with the fact that $|G|\geq \mu$, 
\begin{equation} \label{Sep}
	\min(\hat m_n, 1-\hat m_n)\geq\frac{\mu}{4}.
\end{equation}
Let us define $$I_1^+=\{i\in I_1 : X_i\geq \hat m_n\},$$ and $$I_1^-=\{i\in I_1 : X_i\leq \hat m_n\}.$$ 
By \eqref{Sep}, $\# I_1^\epsilon \geq \frac{\mu n}{8}-1\geq \frac{\mu n}{16}$ for $n$ large enough, and for $\epsilon\in\{+,-\}$. 
Note that $\{X_i : i\in I_1^+\}$ (resp. $\{X_i : i\in I_1^-\}$) is a deterministic and regular design of the segment $[\hat m_n,1]$ (resp. $[0,\hat m_n]$), of cardinality greater or equal to $\frac{\mu n}{16}$, as we saw just before. Then, since we have both
\begin{equation*}
	Y_i=\mathds 1(X_i\leq b)+\xi_i, \forall i\in I_1^+
\end{equation*}
and
\begin{equation*}
	Y_i=\mathds 1(X_i\geq a)+\xi_i, \forall i\in I_1^-,
\end{equation*}
the change points $a$ and $b$ can be estimated as in Theorem \ref{ThmChangePointUB}, using the subsamples $\{(X_i,Y_i) : i\in I_1^+\}$ and $\{(X_i,Y_i) : i\in I_1^-\}$ respectively, and we get two estimators $\tilde a_n$ and $\tilde b_n$ which satisfy:
\begin{equation*}
	\PP_G\left[|\tilde a_n-a|\geq \frac{16y}{\mu n}, E_x\right]\leq C_0e^{-y/(8\sigma^2)}
\end{equation*}
and
\begin{equation*}
	\PP_G\left[|\tilde b_n-b|\geq \frac{16y}{\mu n}, E_x\right]\leq C_0e^{-y/(8\sigma^2)},
\end{equation*}
for all $y>0$.
Set $\tilde G_n=[\tilde a_n,\tilde b_n]$, on the event $E_x$, and $\tilde G_n=\emptyset$ on its complementary $\bar {E_x}$.
By setting $x=\frac{\mu n}{2\ln n}$, which is the maximal value authorized in this proof,
\begin{align*}
	\PP_G\left[|\tilde G_n\triangle G|\geq \frac{y}{n}\right] & \leq \PP_G\left[|\tilde G_n\triangle G|\geq \frac{y}{n}, E_x\right]+\PP_G[\bar{E_x}] \\
	& \leq 2C_0e^{-\mu y/(256\sigma^2)}+C_1n^4e^{-C_2\mu n/2},
\end{align*}
which ends the proof of Theorem \ref{theoremMu}. \hfill $\blacksquare$

\bibliography{Biblio}

\end{document}